 \documentclass[a4paper,11pt]{article}
  \usepackage{amssymb,amsbsy,amsmath,latexsym,theorem,eepic,enumerate,array}
  \usepackage{xypic,times}
  \xyoption{curve}
  \usepackage{geometry,color}
  \def\FTop{\mathsf{FTop}}

\def\cU{\mathcal{U}}
\def\geq{\geqslant}

\def\ecoll{\searrow \hspace{-0.75em}^e\hspace{0.75em}}
\newcommand{\labto}[1]{\stackrel{#1}{\longrightarrow}}
\def\Ce{\mbox{\v{C}}}

\def\<{\langle}
\def\>{\rangle}

\newcommand{\threeaxes}[3]{\def\objectstyle{\scriptstyle}  \objectmargin={0pt}
\xy
(0,0)*+{}="a",(0,-6)*+{\rule{0em}{1.5ex}#2}="b",(7,0)*+{\;#1}="c",
(14,-3)*+{\;#3}="d" \ar@{->} "a";"b" \ar @{->}"a";"c"  \ar
@{->}"a";"d"\endxy }

\newcommand{\directs}[2]{\def\objectstyle{\scriptstyle}  \objectmargin={0pt}
\xy
(0,4)*+{}="a",(0,-2)*+{\rule{0em}{1.5ex}#2}="b",(7,4)*+{\;#1}="c"
\ar@{->} "a";"b" \ar @{->}"a";"c" \endxy }

\newcommand{\xdirects}[2]{\def\objectstyle{\scriptstyle}  \objectmargin={0pt}
\xy
(0,0)*+{}="a",(0,-6)*+{\rule{0em}{1.5ex}#2}="b",(7,0)*+{\;#1}="c"
\ar@{->} "a";"b" \ar @{->}"a";"c" \endxy }
\newcommand{\sdirects}[2]{\def\objectstyle{\scriptstyle}  \objectmargin={0pt}
\xy
(0,2.2)*+{}="a",(0,-2.5)*+{\rule{0em}{1.5ex}#2}="b",(7,2.2)*+{\;#1}="c"
\ar@{->} "a";"b" \ar @{->}"a";"c" \endxy }

\def\rho{\varrho}

\def\A{\alpha}

\def\eps{\varepsilon}
\def\epsilon{\varepsilon}
\def\pt{\partial}

\def\le{\leqslant}
\def\ge{\geqslant}

\def\A{\alpha}

\def\eps{\varepsilon}
\def\epsilon{\varepsilon}
\def\pt{\partial}

\def\le{\leqslant}
\def\ge{\geqslant}

\def\geq{\geqslant}
\def\subset{\subseteq}
\def\ecoll{\searrow \hspace{-0.75em}^e\hspace{0.75em}}

 \def\Im{\mathop{\rm Im}\nolimits}

 \def\eps{\varepsilon}

\def\A{\alpha}

\def\om{$\omega$-}

\def\epsilon{\varepsilon}

\def\red{\textcolor{red}}

\def\blue{\textcolor{blue}}

\raggedright
\def\red{\textcolor{red}}

\def\blue{\textcolor{blue}}
\def\magenta{\textcolor{magenta}}
\def\Ker{\operatorname{Ker}}

\def\io{^{-1}}
  
\begin{document}

  \sffamily \huge

\begin{center}
\red{ \bf \Huge Towards Non Commutative \\ Algebraic Topology}

Ronnie Brown

University of Wales, Bangor
\end{center}
This is slightly edited version of the transparencies for a
seminar at University College London, May 7, 2003. (Not all were
used!)

 References: http://www.bangor.ac.uk/$\sim$mas010

http://www.bangor.ac.uk/$\sim$mas010/fields-art3.pdf

\newpage
 Acknowledgements
to: work of Henry Whitehead; \\ many collaborators, particularly
\\ Philip Higgins, Jean-Louis Loday, Tim
Porter, Chris~Wensley;\\ 21 Bangor research students; \\
Alexander Grothendieck, for correspondence {\it \`a baton rompu}
1982-93 and \red{Pursuing Stacks} (1983, 600 pages).

Current support:  Leverhulme Emeritus Fellowship \\`crossed
complexes and homotopy groupoids' to produce a book with now
agreed title:

\vspace{2cm}

Title of this seminar!

\newpage
Start (1965): My attempt to get a form of the Van Kampen Theorem
(VKT) for the fundamental group which would also calculate the
fundamental group of the circle $S^1$. Discovered Higgins' work on
groupoids!

Motivation: \blue{ expository and aesthetic, thinking about
\red{anomalies}}.

Trying to understand the \red{algebraic structures underlying
homotopy theory.}

Homotopy and deformation \red{underly notions of classification}
in many branches of mathematics.

Aim: \blue{explore the situation} \\ (so can not be directed at
other peoples' problems).

NOT mainstream algebraic topology -- \vspace{1cm}

\red{we are  digging a new and additional channel. }

\newpage
Overall plan is (avoiding the green bit) $$
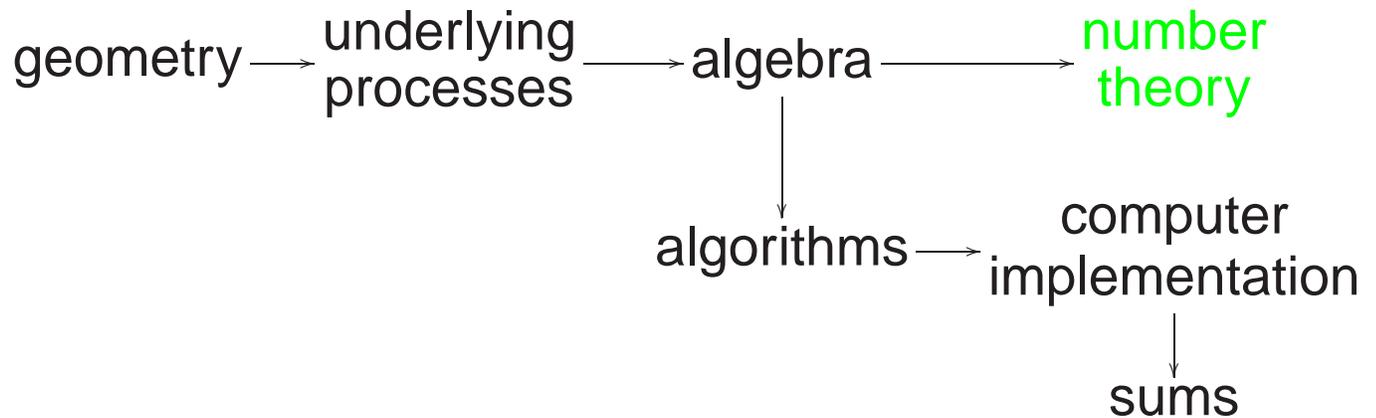
\xymatrix{
\txt{geometry} \ar [r] & \txt{underlying \\ processes} \ar [r] &
\txt{algebra} \ar [d] \ar[r]  &
\textcolor{green}{\txt{number\\theory}} \\& &
 \txt{algorithms} \ar [r] & \txt{computer\\implementation} \ar [d] \\ &&& \txt{sums} }$$

\blue{Example:}  Extensions  $A \rightarrowtail E
\twoheadrightarrow T$ of groups $A,T$ are determined by classes of
\red{factor systems}
$$k^1 : T \to \mbox{Aut}(A), \quad k^2: T \times T \to A.$$
But if $T= gp\langle x,y \mid x^3y^{-2}
\rangle$ is the Trefoil group, then $T$ is infinite, so what can you do? \\
Our result : Extensions are determined by elements $a \in A,\;
a_x,a_y \in \mbox{Aut}(A)$ such that $(a_x)^3(a_y)^{-2}$ is the
inner automorphism determined by $a$.

\blue{Example:} $f: P \to Q$ a morphism of groups. Form the
cofibration sequence $$BP \labto{Bf} BQ \to C(Bf).$$ Find $\pi_2$
and first $k$-invariant of $C(Bf).$ \\ E.g. $(P=C_3 \le Q=S_4)
\implies \pi_2(C(Bf))\cong C_6$. \\
Use \red{non commutative methods} (and computers).

\newpage

Major themes in 20th century mathematics:\\
\begin{tabular}{ll} $\bullet$ \red{non commutativity} & $\bullet$  \red{local to
global}\\  $\bullet$  \red{higher dimensions} & $\bullet$  \red{homology} \\
$\bullet$ \red{$K$-theory} &(Atiyah, Bull LMS, 2002)
\end{tabular}

The VKT for $\pi_1$ is a classical example of a

\red{non commutative local-to-global theorem}.

So any tools which are developed for generalisations and for
higher dimensional forms of it  could be generally useful.

\blue{Non commutative algebraic topology conveniently combines all
the above major themes,} and \red{has yielded substantial }
\magenta{new  calculations, new understanding, new prospects},\\
of which the last is possibly the most important.

Applications to concurrency (GETCO).

Recent EPSRC Grant on \red{Higher dimensional algebra and
differential geometry.}  Peter May's interest in higher
categorical structures. Work with Tony Bak, Tim Porter  on Tony's
`global actions'.

\newpage

Why think of non commutative algebraic topology?

Back in history!

Topologists of the early 20th century knew well that:

1) Non commutative fundamental group $\pi_1(X,a)$ had applications
in geometry, topology,  analysis.

2) Commutative  $H_n(X)$ were defined for all $n \ge 0$.

3) For connected $X$, $$H_1(X) \cong \pi_1(X,a)^{ab}.$$ So they
dreamed of

 \blue{higher dimensional,  non commutative
versions of the fundamental group}.

Gut feeling:  dimensions $> 1$ need invariants which are
\begin{center} \blue{`more non commutative' than groups}.
\end{center}

\newpage

1932:  ICM at Z\"urich:

\v{C}ech:  submits a paper on higher homotopy groups
$$\pi_n(X,x);$$
Alexandroff, Hopf: prove  commutativity for $n \ge 2$;\\
persuade \v{C}ech to withdraw his paper;\\only a small paragraph
appears in the Proceedings.

Reason for commutativity (in modern terms):

\red{group objects in groups are commutative groups.}

$\pi_2(X,x)$, even  considered as a $\pi_1(X,x)$-module, is only a
pale shadow of the 2-type of $X$.

\blue{What is going on? }

\newpage

Overall philosophy: look for \\ \red{ algebraic models of homotopy
types}

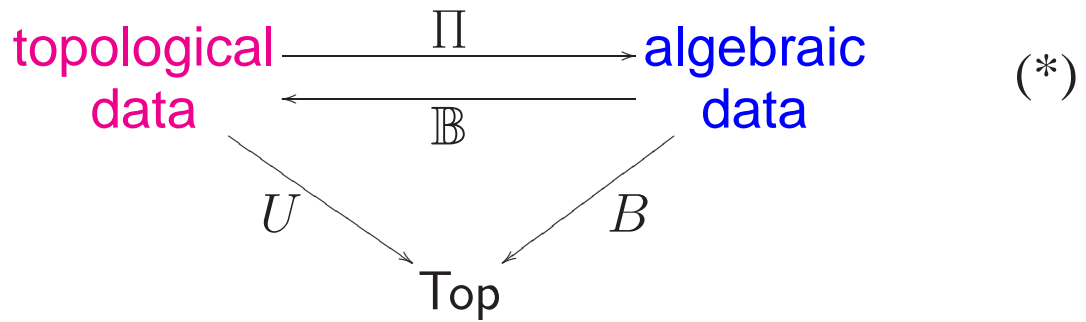
\begin{equation}
\def\labelstyle{\textstyle}
\xymatrix@=4pc{\magenta{\txt{topological\\data}}
\ar [dr] _{U} \ar @<0.75ex>[rr] ^{\Pi} && \ar @<0.75ex> [ll]
^{\mathbb B}\ar [dl] ^B
\blue{\txt{algebraic \\ data}}  \\
&  \mathsf{ Top} &}\tag{*}
\end{equation}

1) $U$ is a forgetful functor and $B= U \circ \mathbb B$;

2) $\Pi$ is defined homotopically;

3) (local to global, allowing calculation!): \\\qquad  $\Pi$
preserves certain colimits;

3) (algebra models the geometry) \\ \qquad $\Pi  \circ{\mathbb B}
\simeq 1$;

4) (capture homotopical information): \\$\exists$ natural
transformation $1 \to{\mathbb B} \circ\Pi$ with good properties.

\newpage
\hspace{5cm}  Some examples:

\begin{tabular}{|m{8cm}|m{7cm}|}\hline
\magenta{topological data} & \blue{algebraic data} \\\hline \hline
\magenta{spaces with base point} & \blue{groups} \\\hline
\magenta{spaces with a set of base points} & \blue{groupoids}
\\ \hline \magenta{based pairs} & \blue{crossed modules} \\ \hline
\magenta{filtered spaces} & \blue{crossed complexes} \\ \hline
\magenta{$n$-cubes of spaces} & \blue{cat$^n$-groups} \\\hline
\end{tabular}

So on the blue side we have various generalisations of groups.
Here are some more!

$$\frame{\def\labelstyle{\textstyle}
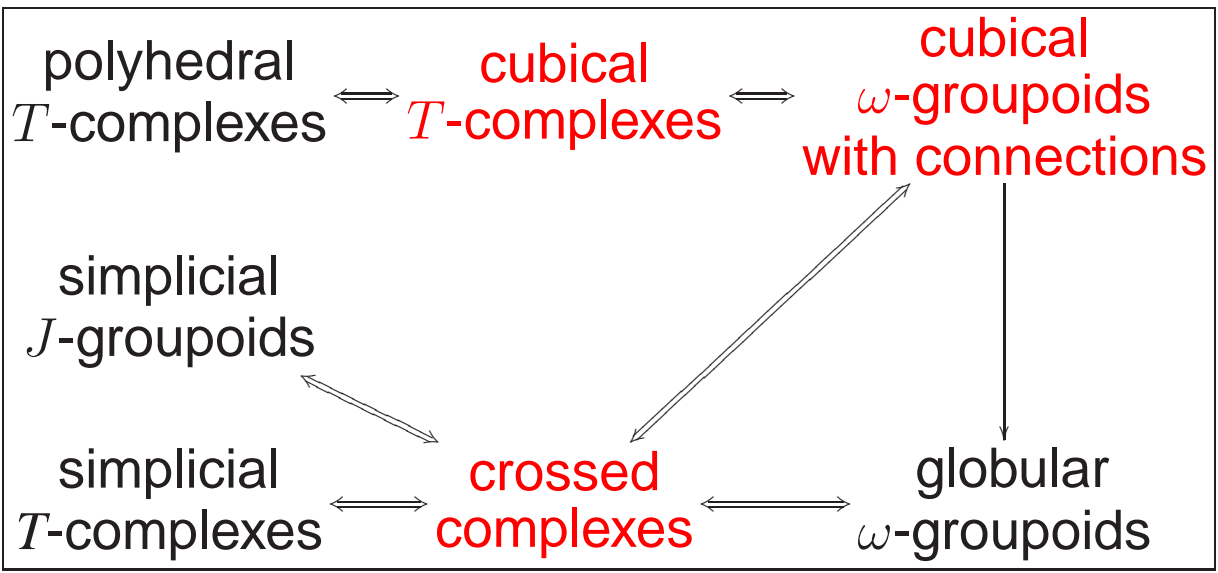
\xymatrix {\txt{polyhedral\\
$T$-complexes} \ar @{<=>} [r] & \red{\txt{cubical
\\$T$-complexes}}   \ar @{<=>} [r] & \red{ \txt{cubical
 \\$\omega$-groupoids\\with connections}} \ar [dd]
\ar @{<=>} [ddl]  \\
 \txt{simplicial \\$J$-groupoids}\ar @{<=>} [dr] && \\
   \txt{simplicial \\{\it  T}-complexes} \ar @{<=>}[r] & \red{ \txt{crossed \\
complexes}} \ar @{<=>}[r] & \txt{ globular \\
$\omega$-groupoids}                                 } }$$

The equivalence of the red structures  is \red{required} for the
proof of the Brown-Higgins GVKT.

\newpage

Features of groupoids: \\ structure in dimension 0 and 1;  \\
 composition operation is partially defined;\\
 allows the combination of groups and graphs, or
groups and space.

\red{Higher dimensional algebra} (for me) is the study and
application of \blue{algebraic structures whose domains of the
operations  are given by geometric  conditions.} This allows for a
vast range of new algebraic structures related to geometry.

Why so many structures?

More compact convex sets in dimension 2 than dimension 1!

The algebra has to express and cope with structures defined by
different geometries.

\newpage

Easiest example: Cubes.

Cubical methods are used in order to express the intuitions of

1) Multiple compositions \\ (algebraic inverses to subdivisions);

2) Defining a commutative cube.

3) Proving a multiple composition of commutative cubes is
commutative (Stokes' Theorem?!).

\blue{4) Construction and properties of higher homotopy
groupoids.}

5) Homotopies and tensor products: $$I^m \times I^n\cong I^{m+n}$$

\newpage
\blue{Category $\FTop$ of filtered spaces:} $$X_*: X_0 \subset X_1
\subset X_2 \subset \cdots \subset X_\infty$$  of subspaces of
$X_\infty$.

\blue{Homotopical quotient:} $$p: R(X_*)\longrightarrow
\rho(X_*)=R(X_*)/ \equiv$$ where
$$R(X_*)_n = \FTop(I^n_*,X_*),$$  $I^n_* = n$-cube with its skeletal
filtration,\\
\red{$\equiv$: homotopy through filtered maps rel vertices of
$I^n$.} \\ Then $R(X_*),\rho(X_*)$ are cubical sets with
connections. \\ (\red{Connections} are extra `degeneracy'
operations.)

But $R(X_*)$ has standard \red{partial compositions:}\\for $i=1,
\ldots, n$, if $a,b \in R(X_*)_n$ and $\partial^+_i a = \partial
^-_i b$ we can define $a \circ_i b\in R(X_*)_n$ by
$$(a \circ_i b)(t_1, \ldots,t_n) = \begin{cases} a(\ldots, 2t_i,
\ldots ) & t_i \le \frac{1}{2}, \\
b(\ldots, 2t_i -1, \ldots ) & t_i \ge \frac{1}{2}.
\end{cases} $$

\newpage
\blue{Major result 1):} The standard compositions  on $R(X_*)$ are
inherited by $\rho(X_*)$ to make it the \red{fundamental cubical
$\omega$-groupoid of $X_*$.}

This is quite difficult to prove, and is non trivial even in
dimension 2. The result is \blue{precise} in that there is just
enough filtration room to prove it.

\blue{Major result 2):} The quotient map $p: R(X_*)\to\rho(X_*)$
is a Kan fibration of cubical sets.

This result is almost unbelievable. Its proof has to give a
systematic method of deforming a cube with the right faces `up to
homotopy' into a cube with exactly the right faces, using the
given homotopies.

\newpage

Here is an application of 2) which is essential in many proofs.

\blue{ Theorem: Lifting multiple compositions }

\red{ Let  $[\alpha_{(r)}]$ be a multiple composition in
$\rho_n(X_*)$. Then representatives  $a_{(r)}$ of the
$\alpha_{(r)}$   may be chosen so that the composition $a_{(r)}$
is well defined in $R_n(X_*)$.}

\blue{Explanation: To say that $[\alpha_{(r)}]$ is well defined
says representatives $a_{(r)}$ agree with neighbours up to
homotopy, and these homotopies are arbitrary. All these homotopies
have to be used to obtain the representatives which actually agree
with their neighbours.}

This is an example of why setting up higher homotopy groupoids is
not straightforward.

\newpage
\blue{ Proof:} The multiple composition $[\alpha_{(r)}]$
determines a cubical map
$$A:K \to\rho(X_*)$$
where the cubical set $K$ corresponds to  a subdivision of the
geometric cube.

Consider the diagram \blue{
$$\xymatrix@=6pc{\ast \ar [r] \ar [d] & R(X_*) \ar [d] ^p \\
K \ar [r] _-A \ar @{-->}[ur]|{A'} & \rho(X_*)}.$$}

Then $K$ collapses to $*$, written $K\searrow *$. \\ By the
fibration result, \\ $A$ lifts to $A'$, which represents
$a_{(r)}$, as required.
\newpage

\blue{Major result 3):} If $X_\infty= U \cup_W V$ with $U,V,W$
open, and the induced filtrations $U_*, V_*, W_*$ are
\red{connected} then

C) $X_*$ is connected;

I) The following diagram

$$\xymatrix{ \rho(W_*) \ar [d] \ar [r] & \rho(V_*) \ar [d] \\
\rho(U_*)\ar [r] & \rho(X_*)}$$ is a pushout of cubical
$\omega$-groupoids with connection.

\blue{Proof Outline:} Verify the universal property with regard to
maps to $G$. Take $ a \in\in  \rho(X_*)_n$. Subdivide
$a=[a_{(r)}]$ so that each $a_{(r)}$ lies in $U$ or $V$. Use
connectivity to deform $a_{(r)}$ to $a_{(r)}'\in R(Y_*), Y= U,V,W$
such that $a'=[a_{(r)}']$ is defined. Map the pieces over to $G$
and recombine. \red{Analogy with email.}

You have to prove independence of choices. This needs a technology
of commutative cubes.

Applications: Translate to crossed complexes.

\newpage

Down to earth and explain \red{crossed modules}

JHC Whitehead in 1939-50 abstracted properties of
\begin{equation}\label{2relhom}  \partial: \pi_2(X,X_1,a) \to
\pi_1(X_1,a) \tag{*} \end{equation} to define a \red{Crossed
Module}:\\  morphism of groups \\  $\mu : M \to P$ and action $M
\times P \to M, \;(m,p)\mapsto m^p$ \\  of the group
$P$ on the group $M$ such that: \\
{\bf CM1)}  $\mu (m^p) = p^{-1}(\mu m)p$ \qquad {\bf CM2)}
$n^{-1}mn= m ^{\mu n}$ \\ for all $m,n \in M, p \in P.$

Now a key concept in  non commutative algebraic topology and
homological algebra.

Simple consequences of the axioms:\\ $\bullet$  $\Im \mu $ is normal in $P$ \\
$\bullet$  $\Ker \mu$ is central in $M$ and is acted on trivially
by $\Im \mu$, so that $\Ker \mu$ inherits an action of $M/\Im
\mu$.
\newpage

Standard algebraic examples:
\par \noindent (i) normal inclusion \red{$ M \triangleleft P$}; \par \noindent (ii) inner
automorphism map  \red{$\chi : M \to \mbox{Aut}\;M $}; \par
\noindent (iii) the zero  map  \red{ $0: M \to P $} where  $M$ is
a $P$-module; \par \noindent (iv) an epimorphism  \red{$M \to P $}
with kernel contained in the centre of $M$.

\noindent {\bf Theorem (Mac~Lane-Whitehead, 1950} {\it Crossed
modules classify all connected weak based homotopy 2-types.}

\blue{Crossed modules as candidates for 2-dimensional groups? }

\newpage

1974 (published 1978): Brown and Higgins proved that the functor
$$\Pi_2: \mbox{(based pairs of spaces)} \to \mbox{(crossed
modules)}$$ preserves certain colimits. This  \blue{allows totally
new 2-dimensional homotopical calculations.} One can compute with
crossed modules in a similar, but more complicated,  manner to
that for groups.

Recent work with Chris Wensley uses symbolic computation to do
more sums.

\blue{The aim of  these new calculations is to \red{prove} (i.e.
test) the power of the machinery.}

Grothendieck's aim in Pursuing Stacks was Non Abelian Homological
Algebra.

The real aim is an \blue{extension of method}, in the belief that
methods last longer than theorems.

Next show examples of a new concept and calculations.

\newpage

\red{Induced crossed modules} (Brown-Higgins, 1978). $f:P \to Q$ a
group morphism.
$$\xy (0,0)*+{\magenta{\txt{crossed \\
$P$-module}}} \endxy \left\{ \vcenter{ \xymatrix@=3pc{M \ar
[d]_\mu \ar [r] &f_*(M) \ar [d] ^\partial
\\ P \ar [r] ^ {f} & Q}} \right\} \xy (0,0)*+{\blue{\txt{crossed \\
$Q$-module}}} \endxy $$

$$ f_*: \mbox{\magenta{crossed $P$-modules}} \to \mbox{\blue{crossed
$Q$-modules}}$$ Universal property: left adjoint to pullback by
$f$.

Construction: generated by symbols \blue{$$m^q, \; m \in M, q \in
Q$$} \noindent  with $\partial (m^q)= q \io (fm) q $ and rules$$
(m^p)^q= m^{(fp)q}, \; \mbox{CM2 for } \partial.$$

Example of a `change of base' construction.
\newpage
Example 1) Let $f:P \to Q$ be a morphism of groups, inducing a
cofibration sequence $$BP \to BQ \to C(Bf).$$ Algebraic
description of the 2-type of $C(Bf)$ as an \red{induced crossed
module} $f_*(P \to P)$, so we can calculate specific examples.

 1) (Brown, Wensley, 1995) $M,P,Q$ finite
$\implies$ $f_*(M)$ finite. Hence \red{computations of homotopy
2-type} of $B(C(Bf))$ when \magenta{$\mu =1_P: P \to P$} and $f: P
\triangleleft Q$; more generally of a homotopy pushout
$$ \xymatrix{BP \ar [d] \ar [r] & BQ \ar [d] \\
B(M \to P) \ar [r] & X }$$

\newpage

2) \magenta{$\mu =1: F(R) \to F(R)$},  $\omega': F(R) \to Q$
defined by $\omega: R \to Q$. Then\blue{
$$\partial: C(\omega)=\omega'_*(F(R)) \to Q$$}\noindent  is the free
crossed $Q$-module on $\omega$. (Defined directly by Whitehead).

Corollary is a major result:

\noindent {\bf Theorem W (1949)}   {\it $$\pi_2(X_1 \cup
\{e^2_r\}_{r \in R},a) \to \pi_1(X_1,a)$$ is isomorphic to the
free crossed $\pi_1(X_1,a)$-module on the classes of the attaching
maps of the 2-cells.
 }

 This is important for relating combinatorial group theory and
 2-dimensional topology. (Identities amon relations. )

 \newpage
\blue{Some Computer Calculations (C.D. Wensley using GAP):}
$[m,n]$ is the $n$th group of order $m$ in GAP.

$M \triangleleft P$;  $f : P \le  S_4$. Calculate $f_\ast M$. \\
Set $C_2 = \<(1,2)\>$, $C_2^{\prime} = \<(1,2)(3,4)\>$, $C_2^2 =
\<(1,2),(3,4)\>$.

\begin{center}

\begin{tabular}{||c|c|c|c|c||}
\hline
$M$ & $P$ & $f_\ast M$ & $\ker \partial$ & $\mbox{Aut}(f_\ast M)$  \\
\hline\hline
$C_2$  & $C_2$       & $GL(2,3)$   & $C_2$   & $S_4 C_2$          \\
$C_3$  & $C_3$     & $C_3\ SL(2,3)$& $C_6$   & $[144,183]$        \\
$C_3$  & $S_3$       & $SL(2,3)$   & $C_2$   & $S_4$              \\
$S_3$  & $S_3$       & $GL(2,3)$   & $C_2$   & $S_4 C_2$          \\
$C_2^{\prime}$ & $C_2^{\prime}$  & $[128,?]$ & $C_4C_2^3$  &      \\
$C_2^{\prime}$ & $C_2^2,C_4$ & $H_8^+$     & $C_4$   & $S_4 C_2$  \\
$C_2^{\prime}$ & $D_8$       & $C_2^3$     & $C_2$   & $SL(3,2)$  \\
$C_2^2$& $C_2^2$     & $S_4C_2$    & $C_2$   & $S_4 C_2$          \\
$C_2^2$& $D_8$       & $S_4$       & $I$     & $S_4$              \\
$C_4$  & $C_4$       & $[96,219]$    & $C_4$   & $[96,227]$        \\
$C_4$  & $D_8$       & $S_4$       & $I$     & $S_4$              \\
$D_8$  & $D_8$       & $S_4 C_2$   & $C_2$   & $S_4 C_2$          \\
\hline
\end{tabular}\\
\end{center}

$\ker \partial \cong \pi_2(C(Bf))$. \\
\red{Need the non commutative structure to find this.} \\ Hard to
determine the first $k$-invariant in $H^3(\mbox{Coker}
\partial,\ker \partial)$.

\blue{Geometric significance of the table?}

\newpage

\blue{Conclusion:}

Key inputs: VKT for the \\ \red{fundamental groupoid
$\pi_1(X,X_0)$ on a set $X_0$ of base points} (RB: 1967).

CLAIM: all of 1-dimensional homotopy theory is better presented
using groupoids rather than groups.

Substantiated in books by Brown (1968) and Higgins (1971). Ignored
by most topologists!

Hint as to higher dimensional prospects:

 \blue{(Group objects in
groupoids) $\Leftrightarrow$ (crossed modules).}
\\ (Grothendieck school, 1960s).
\\ Generalising:\\ \hspace{1em}  (congruences on a group) $\Leftrightarrow$ (normal
subgroups).

\red{Further outlook:} Generalise this to other algebraic
structures than groups.

See work of Fr\"ohlich, Lue, Tim Porter.

Groupoids in  Galois Theories (Grothendieck, Magid, Janelidze).
\newpage

So look for \red{higher homotopy groupoids}.

And applications of groupoids, multiple groupoids, and higher
categorical structures  in mathematics and science.

Hence the term \blue{higher dimensional algebra} (RB, 1987). Web
search shows many applications.

Pursuing Stacks has been a strong international influence.

I gave an invited talk in Delhi in February to an International
Conference on Theoretical Neurobiology!

It is still early days!

\end{document}